\documentclass[12pt]{article}

\usepackage{amsfonts}

\usepackage{amsmath,latexsym,amsthm,amsxtra,amssymb,mathrsfs}
\usepackage{graphics}
\usepackage{enumerate}

\newtheorem{theorem}{Theorem}
\newtheorem{proposition}[theorem]{Proposition}
\newtheorem{lemma}[theorem]{Lemma}

\newtheorem{remark}[theorem]{Remark}

\allowdisplaybreaks

\begin{document}
\begin{center}
\textbf{ A Bando-Mabuchi Uniqueness Theorem}
\end{center}

\centerline{Li YI}
\vskip 0.1cm

\centerline{Institut Elie Cartan, Nancy}
\vskip 0.2cm

\noindent \textbf{Abstract.} \textit{In this paper we prove a uniqueness theorem on generalized K\"ahler-Einstein metrics on Fano manifolds. Our result generalize the one shown by Berndtsson using the convexity properties of Bergman kernels. The same technics as well as that of the regularization of closed positive currents and an estimate due to Demailly-Koll\'ar will be used in our proof.}\\

\noindent \textbf{\S0. Introduction.}\\
\medskip 

\noindent Let $X$ be an $n$-dimensional projective manifold, here $n\geq3;$ we denote by $c_1(X)$ its first Chern class, i.e.
the Chern class corresponding to the anti canonical bundle $-K_X$. Let $\{\alpha\}$ be a 
$(1, 1)$ class on $X$; here we denote by $\alpha$ a non-singular representative. We assume that the 
following requirements are satisfied.

\begin{enumerate} 

\item[$\rm(i)$] The class $c_1(X)- \{\alpha\}$ is K\"ahler, i.e. it contains a K\"ahler metric which we denote by
$\omega$;
\smallskip
\item[$\rm(ii)$] The class $\{\alpha\}$ is pseudo-effective; let $\Theta\in \{\alpha\}$ be a closed positive current; 
\smallskip
\item[$\rm(iii)$] If we denote $\displaystyle {\rm Ric}(\omega)$ the Ricci curvature associated to the 
metric $\omega$, then we have
$${\rm Ric}({\omega})= \omega+ \Theta- 
\sqrt{-1}\partial\overline{\partial}f_\Theta$$
for some function $f_\Theta$ which is unique up to an 
additive constant. We assume that there exists a real number $\varepsilon_0> 0$ such that
$$\int_Xe^{-2(1+ \varepsilon_0)f_\Theta}dV< \infty.$$
In other words, the singularities of the current $\Theta$ are not too wide.

\end{enumerate}

\noindent Following Bedford and Taylor (c.f. \cite{BT}), if $\phi$ is a psh \emph{bounded} function defined on an open subset 
$\Omega\subset {\mathbb C}^n$, then the quantity 
$$\big(\sqrt{-1}\partial\overline{\partial}\phi \big)^n$$
is a positive measure of locally finite mass on $\Omega$. Hence it makes sense to consider the
equation
$$(\omega+ \sqrt{-1}\partial\overline{\partial}\varphi)^n= e^{-\varphi- f_\Theta}\omega^n\leqno (\star)$$
where the solution $\varphi$ is assumed to be bounded, and such that
$$\omega_\varphi:= \omega+ \sqrt{-1}\partial\overline{\partial}\varphi$$ 
is a K\"ahler current (i.e. it is greater than a K\"ahler metric).

We remark that if $\Theta$ is non-singular, then the equation $(\star)$ is equivalent with the
following identity
$${\rm Ric}({\omega_\varphi})= \omega_\varphi+ \Theta.$$
The uniqueness of the solution $\varphi$ up to a biholomorphism of $X$ was established by 
Bando-Mabuchi in \cite{BM}.
If $\Theta$ is the current of integration on a divisor $D$, such that the pair $(X, D)$ is klt,
then the question was settled by Berndtsson in \cite{Ber2}, who found a proof within the framework of the 
psh variation of Bergman kernels.\\

\noindent In this article we show that the following generalization of the result in \cite{Ber2} holds.

\begin{theorem}\label{th:1}
Let $\varphi_0$ and $\varphi_1$ be two bounded solutions of the equation \rm{($\star$)}. Then there exists a biholomorphism $F:X\rightarrow X$ such that 
$$F^*(\omega_{\varphi_1})=\omega_{\varphi_0}\quad\text{and}\quad F^*(\Theta)=\Theta.$$
\end{theorem}
\medskip

If we assume that the anti-canonical class has further positivity properties, for example, $-K_X\geq0,$ then one proves in \cite{Ber2} that the proceeding result holds under the assumption $e^{-f_\Theta}\in L_\text{loc}^1(X).$ However, we remark that in order to have this result, we have to assume two strong positivity properties: $c_1(X)-\{\alpha\}$ is a K\"ahler class, and $c_1(X)$ is a semi-positive class.\\

The arguments we will invoke in order to prove this result generalize the approach by Berndtsson in \cite{Ber2}. It is based on the convexity properties of Bergman kernel, which we recall next.\\

Let $L\rightarrow X$ be a line bundle, which is endowed with a family of non-singular metrics $h_t=e^{-\varphi_t},$ $t\in\mathbb{D},$ the unit disk. We assume that the function 
$$\varphi(t,x)=\varphi_t(x)$$
is non-singular and psh. We can see $L$ as a bundle over the family $X\times\mathbb{D};$ the bundle structure is independent of $t\in\mathbb{D},$ but the hermitian structure varies with respect to $t.$ Let 
$$E:=p_*(K_{X\times\mathbb{D}/\mathbb{D}}+L)$$
be the direct image of the relative adjoint bundle of $L,$ and let $(u_t)_{t\in\mathbb{D}}$ be a holomorphic section of $E.$ Of course, in general $(u_t)$ does not correspond to a holomorphic $n$-form on $X\times\mathbb{D}:$ we can only represent $(u_t)_{t\in\mathbb{D}}$ by an $n$-form $U$ on $X\times\mathbb{D}$ such that $\overline{\partial}U$ is equal to some multiple of $dt,$ and such that $U_{|X\times\{t\}}=u_t.$ In order to evaluate the curvature $\Theta_E$ of $E$ in the direction of the section $(u_t),$ we introduce a family $(v_t)$ of $(n-1)$-forms with values in $L$ such that
\begin{enumerate} [$*$]
\item $\overline{\partial}v_t\wedge\omega=0,$
\item $\partial^{\varphi_t}v_t=P_\perp(\dot{\varphi_t}u_t)$
\end{enumerate} 
where $P$ is the projection on the orthogonal complement of the holomorphic $n$-forms on each fiber, then we have the following formula
$$\langle\Theta_E u_t, u_t\rangle=p_*(c_n\sqrt{-1}\partial\overline{\partial}\varphi_t\wedge\widehat{u}\wedge\overline{\widehat{u}}e^{-\varphi_t})+\int_X\|\overline{\partial}v_t\|^2e^{-\varphi_t}\sqrt{-1}dt\wedge d\overline t$$
where we denote by $\widehat{u}=(u_t-dt\wedge v_t)_t,$ another representation of our initial section $(u_t).$\\

We assume now that for some geometric reasons, the curvature of $E$ vanishes. As a consequence, $v_t$ will be holomorphic, as well as the first term on the right hand side of the formula above. But then we can define a vector field by the formula
$$-v_t=V_t\rfloor u_t;$$
it will be holomorphic outside $u_t=0.$ Actually, under the condition $\Theta_E=0$ one can see that $V_t$ is the complex gradient of $\sqrt{-1}\overline{\partial}\dot{\varphi_t};$ finally, a simple computation shows that the flow associated to $V_1,$ $F,$ is the holomorphic map we seek.\\ 

While trying to implement this scheme in our setting, we have to face two important difficulties. In our case, the bundle $L$ will be just $-K_X,$ endowed with the following family of metrics. We connect $\varphi_0$ and $\varphi_1$ with a geodesic, (this is known to exist), and this will give us a metric $\varphi_t$ whose curvature
$$\omega_t=\omega+\sqrt{-1}\partial\overline{\partial}_X\varphi_t$$
belongs to the class $c_1(X)-\{\alpha\}.$ Then the current 
$$\omega_t+\Theta$$
with $\Theta\in\{\alpha\}$ is closed, positive, and belongs to $c_1(X)$ for each $t.$ Our main concerns are:
\begin{enumerate} [$*$]
\item $\Theta$ is not necessarily smooth;
\item the current $\omega_t=\omega+\sqrt{-1}\partial\overline{\partial}_X\varphi_t$ is not necessarily K\"ahler, it is only known to be positive.
\end{enumerate}
The plan of our article will be as follows: we approximate $\omega_t$ and $\Theta$ with non-singular objects (section 1), and then we construct the analogue $v_t^\nu$ of the $(n-1)$-form $v_t$ mentioned above. Thanks to the fact that $\varphi_0$ and $\varphi_1$ are solutions of $\rm(\star),$ we show that the curvature of $E,$ endowed with the regularized metric, converges to zero. An uniformity argument (which is quite subtle, given that the Lelong set of $\Theta$ could be quite complicated) shows that $v_t^\nu$ converges to a holomorphic $(n-1)$-form. Then the vector field will be defined as above. 
\bigskip

\noindent\textbf{\S1. Regularization process.}\\
\medskip

\noindent Fix a non-singular real $(1,1)$-form $\omega$ in the K\"ahler class $c_1(X)-\{\alpha\},$ a quasi-psh function $\varphi$ is said to be $\omega$-psh if $\omega+\sqrt{-1}\partial\overline{\partial}\varphi\geq0.$ We denote $\rm{PSH}(X,\omega)\cap L^\infty$ the convex space of all $\omega$-psh functions which are bounded on $X.$ Let $\varphi_0$ and $\varphi_1\in\rm{PSH}(X,\omega)\cap L^\infty,$ we call the path of semi-positive forms
$$\omega_t=\omega+\sqrt{-1}\partial\overline{\partial}\varphi_t\quad\text{with}\, t\in[0,1]$$
a (generalized) \textit{geodesic} in $\rm{PSH}(X,\omega).$ The following theorem implies the existence of a continuous geodesic which connects $\omega_0$ and $\omega_1$ (cf. \cite{BD}):
\begin{theorem}\label{th:2}
Assume that the semi-positive closed $(1,1)$-forms $\omega_0$ and $\omega_1$ belong to the same K\"ahler class $c_1(X)-\{\alpha\}$ and have bounded coefficients. Then there exists a geodesic $\omega_t$ connecting $\omega_0$ and $\omega_1$ with the properties that it is continuous on $[0,1]\times X$ and that there is a constant $C$ such that $\omega_t\leq C\omega,$ i.e. $\omega_t$ has uniformly bounded coefficients.
\end{theorem}  

\begin{remark}
The geodesic $\omega_t$ can be constructed in such a way that it is Lipschitz in $t.$ Indeed, like shown in \cite{Ber2}, define first
$$\quad \varphi_t=\sup\{\kappa_t\}\leqno(*)$$
where the supremum is taken over all psh $\kappa_t$ with 
$$\lim_{t\rightarrow 0,1}\kappa_t\leq\varphi_{0,1};$$
then construct a barrier 
$$\chi_t=\text{max}\,(\varphi_0-A\text{Re}t,\varphi_1+A(\text{Re}t-1))$$
with $A$ sufficiently large. This $\chi_t$ is psh and we have
$$\lim_{t\rightarrow 0,1}\chi_t=\varphi_{0,1}.$$
Therefore the supremum in $(*)$ has the same property if we restrict to those $\kappa$ which are larger that $\chi.$ For such $\kappa$ the onesided derivative at $0$ is larger than $-A$ and the onesided derivative at $1$ is smaller than $A.$ We assume from now on that $\kappa$ is independent of the imaginary part of $t,$ $\kappa$ is convex in $t$ so the derivative with respect to $t$ increases, and must lie between $-A$ and $A.$ Therefore $\varphi_t$ satisfies
$$\varphi_0-A\text{Re}t\leq\varphi_t\leq\varphi_0+A\text{Re}t$$
at $0$ and a similar estimate at $1.$ Thus
$$\lim_{t\rightarrow 0,1}\varphi_t=\varphi_{0,1} $$
uniformly on $X.$ In addition, the upper semicontinuous regularization $\varphi^*_t$ of $\varphi_t$ must satisfy the same estimate. Since $\varphi^*_t$ is psh it belongs to the class of competitors for $\varphi_t$ and must coincide with $\varphi_t,$ so $\varphi_t$ is psh.

\end{remark}

\noindent Our proof begins with the regularization for $\varphi_t,$ the ampleness of $c_1(X)-\{\alpha\}$ allows us to write $\varphi_t$ as a uniform limit of a decreasing sequence of non-singular functions $\varphi_t^\nu$ with the property that 
$$\omega_t^\nu:=\omega+\sqrt{-1}\partial\overline{\partial}\varphi_t^\nu\geq0.$$
\medskip

\noindent In order to regularize the current
$$\Theta=\alpha+\sqrt{-1}\partial\overline{\partial}\theta,$$
we use the following result:

\begin{theorem}\rm(\cite{Dem1})\label{th:3}
Let $T$ be a closed almost positive $(1,1)$-current and let $\alpha$ be a non-singular real $(1,1)$-form in the same $\partial\overline{\partial}$-cohomology class as $\Theta,$ i.e. $\Theta=\alpha+\sqrt{-1}\partial\overline{\partial}\theta$ where $\theta$ is an almost psh function. Let $\gamma$ be a continuous real $(1,1)$-form such that $\Theta\geq\gamma.$ Suppose that $\mathcal{O}_{T_X}(1)$ is equipped with a non-singular hermitian metric such that the Chern curvature form satisfies
$$c\left(\mathcal{O}_{T_X}(1)\right)+\pi_{X}^*u\geq 0$$
with $\pi_X:P(T^*X)\rightarrow X$ and with some nonnegative non-singular $(1,1)$-form $u$ on $X.$ Fix a hermitian metric $\omega^0$ on $X.$ Then for every $c>0,$ there is a sequence of closed almost positive $(1,1)$-current $\Theta_{c,k}=\alpha+\sqrt{-1}\partial\overline{\partial}\theta_{c,k}$ such that $\theta_{c,k}$ is non-singular on $X\backslash E_c(\Theta)$ and decreases to $\theta$ as $k$ tends to $+\infty$ (in particular, the current $\Theta_{c,k}$ is non-singular on $X\backslash E_c(\Theta)$ and converges weakly to $\Theta$ on $X$) and such that 
\begin{enumerate} [$(1)$]
\item $\Theta_{c,k}\geq \gamma-\min\{\lambda_k, c\}u-\varepsilon_k\omega^0$ where:
\item $\lambda_k(x)$ is a decreasing sequence of continuous functions on $X$ such that $\lim\limits_{k\rightarrow+\infty}\lambda_k(x)=\nu(\Theta,X)$ at every point.
\item $\varepsilon_k$ is positive decreasing and $\lim\limits_{k\rightarrow+\infty}\varepsilon_k=0.$
\item $\nu(\Theta_{c,k},x)=(\nu(\Theta,x)-c)_+$ at every point $x\in X.$
\end{enumerate}

\end{theorem}
\medskip
\noindent Recall that 
$$E_c(\Theta)=\{x\in X; \nu(\Theta,x)\geq c\}$$
is the $c$-upperlevel set of Lelong numbers of $\Theta.$ A well-known theorem of Siu \cite{Siu} asserts that $E_c(\Theta)$ is an analytic subset of $X$ for any $c>0.$\\

\noindent Choose $u=A\omega^0$ where $A$ is a positive constant large enough and choose $c=\frac{1}{\nu}$ with $\nu\in\mathbb{N}$ a decreasing sequence, according to the above theorem, we get a sequence of functions $\theta^\nu$ such that they are non-singular on $X\backslash E_{\frac{1}{\nu}}(\Theta)$ and
$$\Theta^\nu:=\alpha+\sqrt{-1}\partial\overline{\partial}\theta^\nu\geq-A\min\{\lambda_\nu,\frac{1}{\nu}\}\omega^0-\varepsilon_\nu\omega^0$$
holds on $X$ for every $\nu,$ where $\lambda_\nu(x)$ is a decreasing sequence of continuous functions on $X$ such that $\lim\limits_{\nu\rightarrow+\infty} \lambda_\nu(x)=\nu(\Theta,x)$ at every point and $\varepsilon_\nu$ are positive decreasing such that $\lim\limits_{\nu\rightarrow+\infty}\varepsilon_\nu=0.$\\
\medskip

\noindent During the proof we need to use the "metric" form $"e^{-\cdot}"$ concerning the currents $\omega_t$ and $\Theta.$ So we revise next the notations. Denote $\phi_t=\phi^0+\varphi_t$ with $\sqrt{-1}\partial\overline{\partial}\phi^0=\omega.$
So $\omega_t$ can be written as 
$$\omega_t=\sqrt{-1}\partial\overline{\partial}\phi_t.$$
Similarly, denote $\psi=\psi^0+\theta$ with $\sqrt{-1}\partial\overline{\partial}\psi^0=\alpha.$ Then $\Theta$ can be presented as 
$$\Theta=\sqrt{-1}\partial\overline{\partial}\psi.$$
Hence if we denote $\tau_t=\phi_t+\psi$ and respectively $\tau_t^\nu=\phi_t^\nu+\psi^\nu,$
the above discussions imply that the sequence of currents 
$$\omega_t^\nu+\Theta^\nu=\sqrt{-1}\partial\overline{\partial}\tau_t^\nu$$
decreases to 
$$\omega_t+\Theta=\sqrt{-1}\partial\overline{\partial}\tau_t$$
as $\nu$ tends to $+\infty$ and is non-singular on $X\backslash E_\frac{1}{\nu}(\Theta).$ We evaluate next the lower bounds of $\omega_t^\nu+\Theta^\nu.$ To this end, denote first a neighborhood set of $E_\frac{1}{\nu}(\Theta)$ by
$$U_{\delta_\nu}=\{x\in X: e^{\psi^\nu}<\delta_\nu\}$$
(the choice of $\delta_\nu$ to be precise in a moment). Remark that according to Theorem \ref{th:2}, $\omega_t$ has uniformly bounded coefficients, since $\omega_t^\nu$ converge uniformly to $\omega_t,$ we have actually $e^{\tau_t^\nu}< C\delta_\nu$ for $x\in U_{\delta_\nu}$ with $C>0$ a numerical constant. Then by Theorem \ref{th:3}, we have 
\begin{equation}\label{eq:1}
\sqrt{-1}\partial\overline{\partial}\tau_t^\nu\geq-\left(A\min\{\lambda_\nu,\frac{1}{\nu}\}+\varepsilon_\nu\right)\omega^0
\end{equation}
on $\displaystyle X\backslash U_{\delta_\nu}.$ As $\nu\rightarrow+\infty,$ the quantity $\displaystyle-\left(A\min\{\lambda_\nu,\frac{1}{\nu}\}+\varepsilon_\nu\right)\omega^0$ tends to zero; Globally on $X,$ we have  
\begin{equation}\label{eq:2}
\sqrt{-1}\partial\overline{\partial}\tau_t^\nu\geq-\left(2A\frac{1}{\nu}\right)\omega^0\geq-3A\omega^0.
\end{equation}
\medskip

\noindent In the following proof we need to solve some $\partial^{\tau_t^\nu}$-equations, to insure that the solutions of such equations we will deal with satisfy certain $L^2$-properties, we need to modify further our approximations in order that the new "metrics" are non-singular on $X.$ To this end, define  
$$\widetilde{\tau}_t^\nu=\log\left(e^{\tau_t^\nu}+\frac{1}{\nu}\right).$$
Since $\displaystyle e^{\tau_t^\nu}+\frac{1}{\nu}>0,$ functions $\widetilde{\tau}_t^\nu$ constructed in this way are therefore non-singular on $X.$ We analyze next their Hessian properties.
\begin{align}\label{eq:3}
&\sqrt{-1}\partial\overline{\partial}\log\left(e^{\tau_t^\nu}+\frac{1}{\nu}\right)\nonumber\\
=&\sqrt{-1}\partial\frac{e^{\tau_t^\nu}\overline{\partial}\tau_t^\nu}{e^{\tau_t^\nu}+\frac{1}{\nu}}\nonumber\\
=&\sqrt{-1}\frac{(e^{\tau_t^\nu}\partial\overline{\partial}\tau_t^\nu+e^{\tau_t^\nu}\partial\tau_t^\nu\wedge\overline{\partial}\tau_t^\nu)(e^{\tau_t^\nu}+\frac{1}{\nu})-e^{2\tau_t^\nu}\partial\tau_t^\nu\wedge\overline{\partial}\tau_t^\nu}{(e^{\tau_t^\nu}+\frac{1}{\nu})^2}\nonumber\\
=&\sqrt{-1}\frac{(e^{\tau_t^\nu}+\frac{1}{\nu})e^{\tau_t^\nu}\partial\overline{\partial}\tau_t^\nu+\frac{1}{\nu}e^{\tau_t^\nu}\partial\tau_t^\nu\wedge\overline{\partial}\tau_t^\nu}{(e^{\tau_t^\nu}+\frac{1}{\nu})^2}\nonumber\\
=&\sqrt{-1}\frac{e^{\tau_t^\nu}}{e^{\tau_t^\nu}+\frac{1}{\nu}}\cdot\partial\overline{\partial}\tau_t^\nu+\sqrt{-1}\frac{\frac{1}{\nu}e^{\tau_t^\nu}}{(e^{\tau_t^\nu}+\frac{1}{\nu})^2}\cdot\partial\tau_t^\nu\wedge\overline{\partial}\tau_t^\nu\nonumber.
\end{align}
We hope that Hessians of these new functions $\widetilde{\tau}_t^\nu$ have similar properties as $\tau_t^\nu,$ that is,   
$$\sqrt{-1}\partial\overline{\partial}\widetilde{\tau}_t^\nu$$
would be uniformly bounded from below on $X$ and almost positive (in other words, the lower bounds tend to zero as $\nu$ tends to infinity) on $X\backslash U_{\delta_\nu}.$ Indeed, this will be true if we choose properly the value of $\delta_\nu.$ Observe first that we have 
$$0<\frac{e^{\tau_t^\nu}}{e^{\tau_t^\nu}+\frac{1}{\nu}}\leq1,$$
so the first term in the last equality satisfies our requirements. We now analyze the second term. Notice that we have 
\begin{equation*}
\sqrt{-1}\partial\tau_t^\nu\wedge\overline{\partial}\tau_t^\nu=|\partial\tau_t^\nu|^2\geq0
\end{equation*}
and 
\begin{equation}\label{eq:3}
\frac{\frac{1}{\nu}e^{\tau_t^\nu}}{(e^{\tau_t^\nu}+\frac{1}{\nu})^2}\geq0
\end{equation}
hold on $X.$ So it is sufficient to check that the function in (\ref{eq:3}) goes to zero on $X\backslash U_{\delta_\nu}.$ Indeed, for $x\in X\backslash U_{\delta_\nu},$ we have
$$\frac{\frac{1}{\nu}e^{\tau_t^\nu}}{(e^{\tau_t^\nu}+\frac{1}{\nu})^2}\leq\frac{1}{\nu}\frac{e^{\tau_t^\nu}}{e^{2\tau_t^\nu}}=\frac{1}{\nu}\frac{1}{e^{\tau_t^\nu}} \leq\frac{C}{\nu\delta_\nu}.$$ 
Hence if we choose for example $\displaystyle\delta_\nu=\frac{1}{\sqrt{\nu}},$ we obtain that  
$$\frac{C}{\nu\delta_\nu}=\frac{C}{\sqrt{\nu}}$$
tends to zero as $\nu$ tends to infinity. Therefore the second term fulfills also our requirements. We can conclude from (\ref{eq:1}), (\ref{eq:2}) and the above discussions that 
\begin{equation}\label{eq:4}
\sqrt{-1}\partial\overline{\partial}\widetilde{\tau}_t^\nu\geq-A_1\omega^0
\end{equation}
on $X$ with $A_1>0$ a numerical constant as well as 
\begin{equation}\label{eq:5}
\sqrt{-1}\partial\overline{\partial}\widetilde{\tau}_t^\nu\geq-A_{\delta_\nu,\nu}\omega^0
\end{equation}
on $X\backslash U_{\delta_\nu}$ with the property that $A_{\delta_\nu,\nu}$ tends to zero as $\nu\rightarrow+\infty.$ For the sake of simplifying the notations, we denote from now on $\widetilde{\tau}_t^\nu$ still by $\tau_t^\nu$ with the understanding that it is now a sequence of non-singular functions on $X$ satisfying (\ref{eq:4}) and (\ref{eq:5}) on the relevant sets.

\bigskip 

\noindent\textbf{\S2. Constructions of $v_t^\nu$.}\\
\medskip

\noindent Like explained in the introduction, finding solutions of the $\partial^{\tau_t^\nu}$-equations \begin{equation}\label{eq:6}
\partial^{\tau_t^\nu}v_t^\nu=P^{\nu,t}_\perp(\dot{\tau}_t^\nu u_t)
\end{equation}
is important in the search of holomorphic vector fields. In this expression, $u_t$ is the section of $E$ described in the introduction, 
$$\dot{\tau}_t^\nu=\frac{\partial\tau_t^\nu}{\partial t}$$
and $P_\perp^{\nu,t}$ is the projection on the orthogonal complement of the holomorphic $n$-forms on each fiber,  
and finally the operator $\partial^{\tau_t^\nu}$ is defined by 
$$\partial^{\tau_t^\nu}:=e^{\tau_t^\nu}\partial e^{-\tau_t^\nu}=\partial-\partial\tau_t^\nu\wedge.$$

\noindent We will explain why equation (\ref{eq:6}) can always be solved for $\tau_t^\nu$ fixed in a moment. In order to analyze the uniformity properties of these solutions, the following lemma will be important (it is a slight variation of Lemma 6.2 in \cite{Ber2}).

\begin{lemma}
Let $L$ be a holomorphic line bundle over $X$ with a metric $\xi $ satisfying $\sqrt{-1}\partial\overline{\partial}\xi\geq -c_0\omega$ with $c_0>0$ a numerical constant. Let $\xi_0$ be a smooth metric on $L$ with $\xi\leq\xi_0$ and assume 
$$I:=\int_Xe^{\xi_0-\xi}<\infty.$$
Then there is a constant $B,$ only depending on $I$ and $\xi_0$ such that if $f$ is a $\overline{\partial}$-exact $L$ valued $(n,1)$-form with
\begin{equation}
\int_X|f|^2e^{-\xi} \leq 1,
\end{equation}\label{eq:7}
there is a solution $u$ to $\overline{\partial}u=f$ with
$$\int_X|u|^2e^{-\xi}\leq B.$$
\end{lemma}  
\medskip
\begin{proof}
\noindent The assumptions imply that
$$\|f\|^2_{\xi_0}=\int_X|f|^2e^{-\xi_0}\leq 1.$$
Since $\overline{\partial}$ has closed range for $L^2$-norms defined by smooth metrics, we can solve $\overline{\partial}u=f$ with the estimate 
$$\int_X|u|^2e^{-\xi_0}\leq C_1$$
for some constant $C_1$ depending only on $X$ and $\xi_0.$ This can be shown by contradiction. Indeed, if this is not true, then we have a non-zero sequence of $\overline{\partial}$-exact $L$ valued $(n,1)$-forms $f_n$ with 
$$\|f_n\|^2_{\xi_0}\leq 1,$$
and a sequence of $L$ valued $(n,0)$-forms which solve the $\overline{\partial}$-equation $$\overline{\partial}u_n=f_n$$ such that  
$$\|u_n\|^2_{\xi_0}=n.$$
Since $f_n$ are $\overline{\partial}$-exact and non-zero, those lie $u_n$ are in the space which is orthogonal to the space of $\overline{\partial}$-exact forms. Since the space of $\overline{\partial}$-exact forms is closed, its orthogonal space is also closed. Let us now consider the $(n,0)$-forms
$$\frac{u_n}{\|u_n\|_{\xi_0}},$$
they lie in the space which is orthogonal to $\overline{\partial}$-exact forms, so the limit
\begin{equation}\label{eq:8}
\lim_{n\rightarrow+\infty}\frac{u_n}{\|u_n\|_{\xi_0}}
\end{equation}
must stay in the same space by its closeness. On the other hand, we have 
$$\frac{f_n}{\|u_n\|_{\xi_0}}=\overline{\partial}\frac{u_n}{\|u_n\|_{\xi_0}}$$
with the left hand side goes to zero when $n\rightarrow+\infty.$ In other words, the term in (\ref{eq:8}) is $\overline{\partial}$-exact. That is to say, it lies in the intersection of the two spaces and hence must be equal to zero. But clearly the norm of $\displaystyle\frac{u_n}{\|u_n\|_{\xi_0}}$ is equal to 1, hence its limit cannot be zero and this leads to the contradiction.\\

\smallskip
\noindent Return to the proof of the lemma, in section $\S 1$ we proved that  
$$\sqrt{-1}\partial\overline{\partial}\tau_t^\nu\geq-A_1\omega^0$$
on $X.$ Choose a collection of coordinate balls $B_j$ such that $B_j/2$ cover $X$ and that on the local coordinate sets $z=(z_1,\ldots,z_n)\in B_j,$ the curvature form of $\tau_t^\nu$ satisfies   
$$\sqrt{-1}\partial\overline{\partial}(\tau_t^\nu+A_1|z|^2)\geq 0.$$
By H\"ormander's $L^2$-estimate (cf. \cite{Hor}), on each $B_j$ we can solve $\overline{\partial}u_j=f$ with the estimate 
\begin{equation}\label{eq:9}
\int_{B_j}|u_j|^2e^{-(\tau_t^\nu+A_1|z|^2)} \leq C_2\int_{B_j}|f|^2e^{-(\tau_t^\nu+A_1|z|^2)}\leq C_2
\end{equation}
where $C_2$ depends only on the size of the balls. Here we used again the assumption (\ref{eq:7}) on $B_j$ together with the fact that $|z|^2\leq \text{Vol}\,(B_j)^2.$ Rewrite inequality (\ref{eq:9}), we have
$$\int_{B_j}|u_j|^2e^{-\tau_t^\nu}\cdot e^{-A_1|z|^2}\leq C_2,$$
so on each $B_j,$ this inequality implies
$$\int_{B_j}|u_j|^2e^{-\tau_t^\nu}\leq C_2e^{A_1\rm{Vol}(B_j)^2} :=C_2^\prime$$
with $C_2^\prime$ depends only on the size of the balls. Denote 
$$h_j:=u-u_j,$$
then it is holomorphic on $B_j$ and we have for $\xi_0$ some arbitrary non-singular metric on $L,$ 
$$\int_{B_j}|h_j|^2e^{-\xi_0}\leq C_3.$$
By Cauchy's estimates, this induces in particular that 
$$\sup_{B_j/2}|h_j|^2e^{-\xi_0}\leq C_4.$$
Therefore we get 
$$\int_{B_j/2}|h_j|^2e^{-\tau_t^\nu}\leq\sup_{B_j/2}|h_j|^2e^{-\xi_0}\cdot I\leq C_4 I$$
and finally
\begin{equation}\label{eq:10} 
\int_{B_j/2}|u|^2e^{-\tau_t^\nu}\leq\int_{B_j/2}|h_j|^2e^{-\tau_t^\nu}+\int_{B_j/2}|u_j|^2e^{-\tau_t^\nu}\leq C_5 I,
\end{equation}
where $C_4$ and $C_5$ depend only on the size of balls. Since $B_j/2$ cover $X,$ 
Summing the inequality (\ref{eq:10}) up, we then get uniform estimates for solutions of $\overline{\partial}$-equations, independent of $\nu$ and $t.$
\end{proof}
\medskip

\noindent The reason that we can always solve the $\overline\partial$-equation is that the $\overline\partial$-operator has closed range, this implies that for every $\nu$ and $t,$ the adjoint operator $\overline{\partial}^*_{\tau_t^\nu}$ also has closed range. Since $P_\perp^{\nu,t}(\dot{\tau}_t^\nu u)$ is orthogonal to holomorphic forms, we can solve the equation 
$$\overline{\partial}^*_{\tau_t^\nu}\alpha_t^\nu=P_\perp^{\nu,t}(\dot{\tau}_t^\nu u)$$
with the property that $\overline{\partial}\alpha_t^\nu=0.$ Write $\alpha_t^\nu=v_t^\nu\wedge\omega^1$ where $\omega^1$ is a fixed K\"ahler form on $X,$ then modulo a sign we have 
$$\partial^{\tau_t^\nu}v_t^\nu=\overline{\partial}^*_{\tau_t^\nu}\alpha_t^\nu.$$ 
In other words, $v_t^\nu$ solve the equation (\ref{eq:6}). The fact that the $L^2$-norm of $v_t^\nu$ is uniformly bounded can be concluded from the above lemma together with some functional analysis, as we will show in the next. Indeed, given an $(n,1)$-form $\alpha_t^\nu$ such that $\|\alpha\|_{\tau_t^\nu}<+\infty,$ by the above lemma there exists an $(n,0)$-form $f_t^\nu$ which solves the $\overline{\partial}$-equation 
$$\overline{\partial}f_t^\nu=\alpha_t^\nu$$
together with the estimate 
\begin{equation}\label{eq:11}
\|f_t^\nu\|_{\tau_t^\nu}\leq C\|\alpha_t^\nu\|_{\tau_t^\nu}.
\end{equation}
with $C$ a constant independent of $t$ and $v.$ Denote
$$\langle\langle\alpha_t^\nu,\alpha_t^\nu\rangle\rangle_{\tau_t^\nu}=\int_X|\alpha_t^\nu|^2e^{-\tau_t^\nu},$$
then  
$$\langle\langle\alpha_t^\nu,\alpha_t^\nu\rangle\rangle_{\tau_t^\nu}=\langle\langle\overline{\partial}^*_{\tau_t^\nu}\alpha_t^\nu,f_t^\nu\rangle\rangle_{\tau_t^\nu}\leq C\|\alpha_t^\nu\|_{\tau_t^\nu}\|\overline{\partial}^*_{\tau_t^\nu}\alpha\|_{\tau_t^\nu}.$$
Rewrite $\alpha_t^\nu=v_t^\nu\wedge\omega^1$ with $\omega^1$ a fixed K\"ahler form on $X,$ we obtain that
$$\|v_t^\nu\|_{\tau_t^\nu}\leq C\|\overline{\partial}^*_{\tau_t^\nu}\alpha_t^\nu\|_{\tau_t^\nu}=C\|P_\perp^{t,\nu}(\dot{\tau}_t^\nu u_t)\|_{\tau_t^\nu}\leq C\|\dot{\tau}_t^\nu u_t\|_{\tau_t^\nu}$$
with the same constant $C$ as in (\ref{eq:11}). By the definition of $\tau_t^\nu,$ 
$$\dot{\tau}_t^\nu=\frac{\partial}{\partial t}(\phi_t^\nu+\psi^\nu)=\dot{\phi}_t^\nu$$
Since $\psi^\nu$ is independent of $t.$ By construction, $\phi_t^\nu$ is uniformly bounded with respect to $\nu$ and Lipschitz in $t,$ so $\dot{\phi}_t^\nu$ is uniformly bounded with respect to both $\nu$ and $t$ and $|\tau_t^\nu-\tau_0^\nu|\leq C.$ Hence
$$\|\dot{\tau}_t^\nu u_t\|_{\tau_t^\nu}$$
is uniformly bounded and this asserts that the $L^2$-norm of $v_t^\nu$ is uniformly bounded.  
\bigskip

\noindent\textbf{\S3. Finding holomorphic vector fields.}\\
\medskip

\noindent Apart from the $L^2$ requirements, we need the limit of $v_t^\nu$ to be holomorphic in order to construct the holomorphic vector field. So in the following we are going to show that the sequence of $L^2$-norms of $\overline{\partial}v_t^\nu$ tends to zero on $X.$ To this end, notice that for every $\nu,$ we can establish a formula concerning $v_t^\nu$ thanks to the following theorem:  
\begin{theorem}\rm(\cite{Ber2}) \label{th:5}
Let $\Theta_E$ be a curvature form on $E$ and let $u_t$ be a holomorphic section of $E.$ For each $t$ in $\omega$ let $v_t$ solve
$$\partial^{\tau_t}v_t=\pi_\perp(\dot{\tau}_tu_t)$$
and be such that $\overline{\partial}_Xv_t\wedge\omega=0.$ Put
$$\widehat{u}=u_t-dt\wedge v_t.$$
Then 
$$\langle\Theta_E u_t,u_t\rangle_t=p_*(c_n\sqrt{-1}\partial\overline{\partial}\tau\wedge\widehat{u}\wedge\overline{\widehat{u}}e^{-\tau})+\int_{X}\|\overline{\partial}v_t\|^2e^{-\tau_t}\sqrt{-1}dt\wedge d\overline t.$$\qed
\end{theorem} 

\noindent More precisely, define
$$\mathcal{F}(t):=-\log\int_Xe^{-\tau_t}d\lambda,$$
and respectively 
$$\mathcal{F}_\nu(t):=-\log\int_Xe^{-\tau_t^\nu}d\lambda,$$
we can apply Theorem \ref{th:5} to $u_t=u$ and obtain the following formula 
\begin{equation}\label{eq:12}
\|u\|^2_{t,\nu}\sqrt{-1}\partial\overline{\partial}\mathcal{F}_\nu=\langle\Theta_E^\nu u,u\rangle_t=
p_*\left(c_n\sqrt{-1}\partial\overline{\partial}
\tau_t^{\nu}\wedge\widehat{u}\wedge\overline{\widehat{u}}e^{-\tau_t^\nu}\right)+\int_X\|\overline{\partial}v_t^\nu\|^2e^{-\tau_t^\nu}\sqrt{-1}dt\wedge d\overline t.
\end{equation}

\medskip 
\noindent This formula connects the $L^2$-norm of $\overline{\partial}v_t^\nu$ and $\displaystyle \sqrt{-1}\partial\overline{\partial}\mathcal{F}_\nu(t).$ If $\mathcal{F}(t)$ is affine along the geodesic, then we have  
$$\sqrt{-1}\partial\overline{\partial}\mathcal{F}(t)=0;$$
this implies that 
$$\sqrt{-1}\partial\overline{\partial}\mathcal{F}_\nu(t)$$
goes to zero weakly on $\mathbb{D}.$ It is the starting point to deduce that $\overline{\partial}v_t^\nu$ goes to zero. Before going further, we first show that $\mathcal{F}(t)$ is affine provided that $\varphi_0$ and $\varphi_1$ are two solutions of equation ($\star$). Here we cite the relevant arguments in \cite{Ber2} for the coherency of the paper.\\

\noindent If $\varphi_0$ and $\varphi_1$ are two such local potentials of the currents $\omega_0$ and $\omega_1$ respectively, we connect them by a continuous geodesic 
$$\omega_t=\omega+\sqrt{-1}\partial\overline{\partial}\varphi_t=\sqrt{-1}\partial\overline{\partial}\phi_t.$$
Since $\omega_t$ is semi-positive, the relative energy
$$\mathcal{E}(\phi_0,\phi_1)$$
is well defined. If $\phi_t$ is non-singular in $t,$ further calculations give that
$$\frac{\partial}{\partial t}\mathcal{E}(\phi_t,\phi_1)=\int_X\dot{\phi}_t(\sqrt{-1}\partial\overline{\partial}\phi_t)^n/\|T\|$$
and 
$$\sqrt{-1}\partial\overline{\partial}_t\mathcal{E}(\phi_t,\phi_1)=\frac{p_*((\sqrt{-1}\partial\overline{\partial}_{X,t}\phi)^{n+1})}{\|T\|}.$$
In the above expression, $p$ is the projection map from $X\times\mathbb{D}$ to $\mathbb{D}$ and $\|T\|$ is the normalizing factor 
$$\|T\|=\int_X(\sqrt{-1}\partial\overline{\partial}_X\phi)^n,$$
chosen so that the derivative of $\mathcal{E}$ becomes $1$ if $\phi_t=\phi+t.$ As shown for example in \cite{Ber},\cite{Mab} and \cite{Sem}, 
$$(\sqrt{-1}\partial\overline{\partial}_{X,t}\phi)^{n+1}=c(\phi)\wedge\sqrt{-1}dt\wedge d\overline{t}\wedge(\sqrt{-1}\partial\overline{\partial}\phi)^n$$
with 
$$c(\phi)=\frac{\partial^2\varphi}{\partial t\overline{\partial}t}-\left|\overline{\partial}_X\frac{\partial\varphi}{\partial t}\right|^2_{\sqrt{-1}\partial\overline{\partial}_X\phi}.$$
Therefore we have
$$\sqrt{-1}\partial\overline{\partial}_t\mathcal{E}(\phi_t,\phi_1) =\frac{\sqrt{-1}dt\wedge d\overline{t}\int_Xc(\phi_t)(\sqrt{-1}\partial\overline{\partial}_X\phi_t)^n}{\|T\|}.$$
It is also shown in these references that the function $c(\phi)$ is equal to the geodesic curvature of the path defined by $\phi$ in the space of K\"ahler potentials. So in particular, if $\phi$ solves the homogenous Monge-Amp\`ere equation
$$(\sqrt{-1}\partial\overline{\partial}_{X,t}\phi)^{n+1}=0,$$
or equivalently $c(\phi)=0,$ the function $\mathcal{E}(\phi_t,\phi_1)$ is harmonic in $t.$ Hence this function is linear along geodesics.

For the family $\phi_t$ which is only bounded, the above calculations are done by defining for every $t,$ 
$$\mathcal{E}(\phi_t,\phi_1)=\inf\{\mathcal{E}(\widetilde{\phi_t},\phi_1)|\,\widetilde{\phi_t}\, \text{are}\,\text{non-singular,}\, \widetilde{\phi_t}\geq\phi_t\}.$$

\noindent Define $$\mathcal{G}_\psi(\phi_t)=-\log\int_Xe^{-(\phi_t+\psi)}d\lambda-\mathcal{E}(\phi_t,\chi)=\mathcal{F}(t)-\mathcal{E}(\phi_t,\chi)$$
with $\chi$ arbitrary and $\psi$ fixed. Then $\phi_0$ and $\phi_1$ are critical points of $\mathcal{G}.$ Indeed, the fact that $\phi_t$ depends only on the real part of $t$ implies that $\mathcal{G}_\psi(\phi_t)$ is convex. Now since $\mathcal{F}$ is convex, $\phi$ is convex in $t,$ we see that the onesided derivatives of $\mathcal{F}(t)$ at the endpoints equal
$$\int\dot{\phi}_te^{-\tau}d\lambda/\int e^{-\tau}d\lambda,$$
with $\dot{\phi_t}$ now stands for the onesided derivatives. On the other hand, the function $\mathcal{E}(\phi_t,\chi)$ is linear, so its distributional derivative
$$\int_X\dot{\phi}_t(\sqrt{-1}\partial\overline{\partial}\phi_t)^n/\|T\|$$
is constant, and hence is equal to its values at the end points by applying the simple convergence theorems for the Monge-Amp\`ere operator. Hence both endpoints are critical points and therefore the convexity of $\mathcal{G}_\psi$ implies that it is constant. Finally, since $\mathcal{E}$ is affine along the geodesic, it follows that 
$$t\rightarrow\mathcal{F}(t)$$
is also affine along the geodesic.\\
\medskip

\noindent With this consideration, we now show that $\overline{\partial}_Xv_t^\nu$ goes to zero weakly over $X\times K$ for any compact $K\subset\mathbb{D}.$ Firstly, we know that for every $t,$ $\tau_t^\nu$ decrease to $\tau_t$ and hence
$$\sqrt{-1}\partial\overline{\partial}\mathcal{F}_\nu$$
goes weakly to
$$\sqrt{-1}\partial\overline{\partial}\mathcal{F}=0$$ 
on $\mathbb{D}.$ By curvature formula (\ref{eq:12}), we see that its right hand side terms 
\begin{equation}\label{eq:13}
 p_*(c_n\sqrt{-1}\partial\overline{\partial}
\tau_t^{\nu}\wedge\widehat{u}\wedge\overline{\widehat{u}}e^{-\tau_t^\nu})+\int_X\|\overline{\partial}v_t^\nu\|^2e^{-\tau_t^\nu}\sqrt{-1}dt\wedge d\overline t
\end{equation}
goes to zero weakly as well. We aim to show that the second quantity in (\ref{eq:13}),
\begin{equation}\label{eq:14}
\int_X\|\overline{\partial}v_t^\nu\|^2e^{-\tau_t^\nu}\sqrt{-1}dt\wedge d\overline t
\end{equation}
goes to zero. To this end, let us examine first the first quantity
$$ p_*(c_n\sqrt{-1}\partial\overline{\partial}
\tau_t^{\nu}\wedge\widehat{u}\wedge\overline{\widehat{u}}e^{-\tau_t^\nu}).$$
As shown in $\S 1,$ it is uniformly bounded from below on $X$ by  
$$-A_1\|\widehat{u}\|^2.$$
This implies in particular that the second quantity (\ref{eq:14}) is at least uniformly bounded from above. More explicit estimates give that for any $t,$ the first quantity is bounded from below by 
\begin{equation}\label{eq:15}
-A_{\delta_\nu,\nu}\|\widehat{u}\|^2-A_1\int_{U_{\delta_\nu}}|v_t^\nu|^2e^{-\tau_t^\nu}
\end{equation}
with $A_{\delta_\nu,\nu}$ tends to zero as $\nu\rightarrow+\infty.$ So we need to evaluate the quantity
\begin{equation}\label{eq:16} 
\int_{U_{\delta_\nu}}|v_t^\nu|^2e^{-\tau_t^\nu}
\end{equation}
when $\nu\rightarrow+\infty.$ We remark that the difference from the case shown in \cite{Ber2} is that here the domain of integration $U_{\delta_\nu}$ changes with $\nu.$ Nevertheless, as we discussed in $\S 2,$ the integral 
\begin{equation}\label{eq:17} 
\int_X|v_t^\nu|^2e^{-\tau_t^\nu}
\end{equation}
as well as
$$\int_X|{\partial}v_t^\nu|^2e^{-\tau_t^\nu}$$
are uniformly bounded from above by constants independent of $t$ and $\nu;$ we just mentioned that the quantity
$$\int_X|\overline{\partial}v_t^\nu|^2e^{-\tau_t^\nu}$$
is also uniformly bounded from above. This asserts that 
$$\int_{X}|dv_t^\nu|^2\leq\int_X|dv_t^\nu|^2e^{-\tau_t^\nu}\leq2\int_X|{\partial}v_t^\nu|^2e^{-\tau_t^\nu}+2\int_X|\overline{\partial}v_t^\nu|^2e^{-\tau_t^\nu}$$
is also uniformly bounded. We can then apply the Sobolev-Poincar\'e type inequality on $X$ (c.f. \cite{HV}) to $v_t^\nu$ and obtain that the quantity  
$$\int_X|v_t^\nu|^\frac{2n}{n-2}d\lambda$$
is uniformly bounded. For $n\geq4,$ we have in particular 
$$\int_X|v_t^\nu|^4d\lambda\leq \left(\int_X|v_t^\nu|^\frac{2n}{n-2}d\lambda\right)^\frac{2(n-2)}{n}\leq C.$$   

\noindent H\"older inequality the yields the following:
\begin{equation}\label{eq:18}
\int_{U_{\delta_\nu}} |v_t^\nu|^2e^{-\tau_t^\nu}\leq\nu\int_{U_{\delta_\nu}}|v_t^\nu|^2d\lambda=\nu\int_{X}1_{U_{\delta_\nu}}|v_t^\nu|^2d\lambda\leq\nu\left(\int_X|v_t^\nu|^4d\lambda\right)^\frac{1}{2}\left(\text{Vol}\,(U_{\delta_\nu})\right)^\frac{1}{2}.
\end{equation}

\noindent Recall that 
$$U_{\delta_\nu}=\left\{x\in X: e^{\tau_t^\nu}<\frac{C}{\sqrt{\nu}}\right\}.$$
We will apply the following proposition to evaluate the volume of $U_{\delta_\nu}:$

\begin{proposition} \rm(\cite{DK}) 
Let $X$ be a complex manifold and $\varphi$ a psh function of $x.$ Let $K\subset X$ be a compact set, $U\Subset X$ a relatively compact neighborhood of $K,$ and let $\mu_U$ be the Riemannian measure on $U$ associated with some choice of hermitian metrics $\omega$ on $X.$ Then for any $c\geq 0,$ the volume of sublevel sets $\{\varphi<\log r\}$ satisfies that 
$$r^{-2c}\mu_U(\{\varphi<\log r\})\leq\int_Ue^{-2c\varphi}dV_\omega.$$
\end{proposition}  
\bigskip 

\noindent According to the above proposition, the volume of $U_{\delta_\nu}$ satisfies the following inequality: For $\beta\geq0,$ we have 
$$\text{Vol}\,(U_{\delta_\nu})\leq\left(\frac{C}{\sqrt{\nu}}\right)^{2\beta}\int_Xe^{-2\beta\tau_t^\nu}d\lambda\leq\left(\frac{C}{\sqrt{\nu}}\right)^{2\beta}\int_Xe^{-2\beta\tau_t}d\lambda.$$
So if we choose for example, $\beta<(1+\varepsilon_0),$ then there exists a constant $B>0$ such that
$$\int_Xe^{-2\beta\tau_t} d\lambda\leq B.$$

\noindent Therefore, the quantity (\ref{eq:18}) reads as 
$$\int_{U_{\delta_\nu}}|v_t^\nu|^2e^{-\tau_t^\nu}\leq\frac{C}{\nu^{\varepsilon_0}}$$
who tends to zero when $\nu$ tends to infinity.\\

\noindent For $n=3,$ the Sobolev-Poincar\'e inequality affirms that 
$$\int_X |v_t^\nu|^6$$
is uniformly bounded from above. Apply H\"older inequality to the lhs quantity in (\ref{eq:18}) now gives that 
$$\int_{U_{\delta_\nu}} |v_t^\nu|^2e^{-\tau_t^\nu}\leq\nu\left(\int_X|v_t^\nu|^6\right)^\frac{1}{3}\left(\text{Vol}\,(U_{\delta_\nu})\right)^\frac{2}{3}\leq\frac{C}{\nu^{\frac{4}{3}(1+\varepsilon_0)-1}}.$$

\bigskip

\noindent All in all, as a consequence, (\ref{eq:16}) hence (\ref{eq:15}) tend to zero. We have proved actually that the first quantity in (\ref{eq:13}) is positive when $\nu$ tends to infinity. By the fact that the second quantity in (\ref{eq:13}) is positive and $\sqrt{-1}\partial\overline{\partial}\mathcal{F}_\nu$ tends to zero, we conclude that both quantities in (\ref{eq:13}) tend to zero. In particular, the $L^2$-norm of $\overline{\partial}_Xv_t^\nu$ goes to zero weakly. Since the $L^2$-norm of $v_t^\nu$ is uniformly bounded we can select a subsequence of $v_t^\nu$ that converges weakly to a form $v_t$ in $L^2.$ Therefore we get
$$\overline{\partial}_Xv_t=0$$
and 
$$\partial_X^{\tau_t}v_t=P^t_\perp(\dot{\tau}_tu)$$
in the weak sense that
$$\int_{X\times\mathbb{D}}dt\wedge d\overline{t}\wedge v_t\wedge\overline{\overline{\partial}W}e^{-\tau_t} =\int_{X\times\mathbb{D}}dt\wedge d\overline{t}\wedge P^t_\perp(\dot{\tau_t}u)\wedge\overline{W}e^{-\tau_t}$$  
for any compact supported non-singular form $W$ of appropriate degree. Indeed, by integration by parts, for every $\nu,$ we have 
$$\int_{X\times\mathbb{D}}dt\wedge d\overline{t}\wedge v_t^\nu\wedge\overline{\overline{\partial}W}e^{-\tau_t^\nu}=\int_{X\times\mathbb{D}}dt\wedge d\overline{t}\wedge\partial^{\tau_t^\nu}(v_t^\nu)\wedge\overline{W}e^{-\tau_t^\nu}=\int_{X\times\mathbb{D}}dt\wedge d\overline{t}\wedge P^{t,\nu}_\perp(\dot{\tau_t^\nu}u)\wedge\overline{W}e^{-\tau_t^\nu}.$$
The left hand side integral in the above equalities converges to 
$$\int_{X\times\mathbb{D}}dt\wedge d\overline{t}\wedge v_t\wedge\overline{\overline{\partial}W}e^{-\tau_t}.$$
We now deal with the wright hand side integral, since $\tau_t^\nu$ is uniformly Lipschitz in $t,$ $\dot{\tau_t^\nu}$ is uniformly bounded; since moreover $\tau_t^\nu$ is convex with respect to $t,$ $\dot{\tau_t^\nu}$ increase. The following lemma implies that $\dot{\tau_t^\nu}$ converge to $\dot{\tau_t}.$  
\begin{lemma}\label{le1}
Let $f_\nu$ be a sequence of non-singular convex functions on an interval in $\mathbb{R}$ that converge uniformly to the convex function $f.$ Let $a$ be a point in the interval such that $f^\prime(a)$ exists. Then $f_\nu^\prime(a)$ converge to $f^\prime(a).$ Since $a$ convex function is differentiable almost everywhere it follows that $f^\prime_\nu$ converges to $f^\prime$ almost everywhere, with dominated convergence on any compact subinterval.
\end{lemma} 
\noindent By the above lemma, $\dot{\tau_t^\nu}u$ converge to $\dot{\tau_t}u.$ Since the sequence $\dot{\tau_t^\nu}u$ are increasing, for any $t,$ we have 
$$|P_\perp^{t,\nu}(\dot{\tau_t^\nu}u)|\leq|\dot{\tau_t^\nu}u|\leq|\dot{\tau_t}u|.$$
That is to say, $P_\perp^{t,\nu}(\dot{\tau_t^\nu}u)$ is uniformly bounded. By Banach-Alaoglu-Bourbaki theorem, there exists a subsequence of $P_\perp^{t,\nu}(\dot{\tau_t^\nu}u)$ which converges weakly, say, to $P_\perp^t(\dot{\tau_t}u).$ In other words, there exists a subsequence of $P_\perp^{t,\nu}(\dot{\tau_t^\nu}u)$ such that 
$$\int_{X\times\mathbb{D}}dt\wedge d\overline{t}\wedge P^{t,\nu}_\perp(\dot{\tau_t^\nu}u)\wedge\overline{W}e^{-\tau_t^\nu}\rightarrow\int_{X\times\mathbb{D}}dt\wedge d\overline{t}\wedge P^t_\perp(\dot{\tau_t}u)\wedge\overline{W}e^{-\tau_t}.$$ 
Hence 
$$\int_{X\times\mathbb{D}}dt\wedge d\overline{t}\wedge v_t\wedge\overline{\overline{\partial}W}e^{-\tau_t} =\int_{X\times\mathbb{D}}dt\wedge d\overline{t}\wedge P^t_\perp(\dot{\tau_t}u)\wedge\overline{W}e^{-\tau_t}.$$\\  
\medskip

\noindent The proof ends if there are no nontrivial holomorphic vector fields on $X.$ Then $v$ must be zero, so $\dot{\tau}_t$ is holomorphic, hence constant. Therefore 
$$\partial\overline{\partial}\dot{\tau}_t=0$$
so $\partial\overline{\partial}\tau_t$ does not depend on $t.$ In the general case, we show that $v_t$ is holomorphic in $t.$ Since we don't know any regularity of $v_t$ except that it lies in $L^2,$ we need to formulate holomorphicity weakly. That is, we are going to prove
$$\int_{X\times\mathbb{D}}dt\wedge\overline{\partial}_t v_t\wedge\overline{\gamma}e^{-\tau_t^\prime}=0$$
for any non-singular $\tau_t^\prime$ on $-K_X$ and $\gamma$ an $(n, 1)$ form on $X\times\mathbb{D}$ which does not contain $dt.$ For this, we follow the same steps as shown in section 4 of \cite{Ber2}.\\
\medskip

\noindent Denote
$$a_\nu:=\int_{X\times\mathbb{D}^\prime}\sqrt{-1}\partial\overline{\partial}\tau_t^\nu\wedge\widehat{u}\wedge\overline{\widehat{u}}e^{-\tau_t^\nu} ,$$
we have discussed that it goes to zero if $\mathbb{D}^\prime$ is a relatively compact subdomain in $\mathbb{D}.$ Shrinking $\mathbb{D}$ slightly we assume that this actually holds with $\mathbb{D}^\prime=\mathbb{D}.$ Choose $W$ to contains no differential $dt,$ so it is an $(n,0)$-form on $X$ with coefficients depending on $t,$ then we have  
$$\int_{X\times\mathbb{D}}\sqrt{-1}\partial\overline{\partial}\tau_t^\nu\wedge\widehat{u}\wedge\overline{W}e^{-\tau_t^\nu}=\int_{X\times\mathbb{D}}\sqrt{-1}\partial\overline{\partial}_t\tau_t^\nu\wedge\widehat{u}\wedge\overline{W}e^{-\tau_t^\nu}.$$
Since $$\sqrt{-1}\partial\overline{\partial}_t\tau_t^\nu\geq 0,$$
by Cauchy-Schwarz inequality, we have 
\begin{equation}\label{eq:19}
\left|\int_{X\times\mathbb{D}}\sqrt{-1}\partial\overline{\partial}\tau_t^\nu\wedge\widehat{u}\wedge\overline{W}e^{-\tau_t^\nu}\right|^2\leq a_\nu\int_{X\times\mathbb{D}}\sqrt{-1}\partial\overline{\partial}_t\tau_t^\nu\wedge W\wedge\overline{W}e^{-\tau_t^\nu}.
\end{equation}
Assume now $W$ has compact support. We apply the following one variable H\"ormander's theorem to the right hand side integral in (\ref{eq:19}) to get that it is dominated by
\begin{equation}\label{eq:20}
\int_{X\times\mathbb{D}}|\partial_t^{\tau^\nu}W|^2e^{-\tau_t^\nu}.
\end{equation}
\begin{proposition} (One variable H\"ormander's theorem). Let $\Omega$ be a domain in $\mathbb{C}$ and let $\phi\in C^2(\Omega),$ let $\alpha\in C_c^2(\Omega).$ Then 
$$\int\Delta|\alpha|^2e^{-\phi}+\int|\frac{\partial}{\partial\overline{z}}\alpha|^2e^{-\phi}=\int|\overline{\partial}_\phi^*\alpha|^2e^{-\phi}.$$\qed
\end{proposition}

\noindent Assume moreover that $W$ is Lipschitz with respect to $t$ as a map from $\mathbb{D}$ into $L^2(X).$ Then (\ref{eq:20}) is uniformly bounded by Lemma \ref{le1}. Hence the left hand side integral in (\ref{eq:19}) goes to zero, in particular, the part which contains $dt\wedge d\overline{t}$ 
$$\int_{X\times\mathbb{D}}\sqrt{-1}dt\wedge d\overline{t}\wedge\left(\partial\overline{\partial}_t\tau_t^\nu-\partial_X\left(\frac{\partial\tau_t^\nu}{\partial\overline{t}}\right)(V_t)\right)\overline{W}e^{-\tau_t^\nu}$$
goes to zero, where $V_t$ is a vector field defined by
$$-v_t=V_t\rfloor u.$$
 Denote 
 $$\mu_t^\nu=\partial\overline{\partial}_t\tau_t^\nu-\partial_X\left(\frac{\partial\tau_t^\nu}{\partial\overline{t}}\right)(V_t).$$
 The next lemma shows that 
$$\int_{X\times\mathbb{D}}\sqrt{-1}dt\wedge d\overline{t}\wedge\mu_t^\nu\wedge\overline{P_\perp W}e^{-\tau_t^\nu}$$
also goes to zero.
\begin{lemma}\label{le2}
Let $\alpha_t$ be forms on $X$ with coefficients depending on $t$ on $\mathbb{D}.$ Assume that $\alpha_t$ is Lipschitz with respect to $t$ as a map from $\mathbb{D}$ to $L^2(X).$ Let $P^t$ be the orthogonal projection on $\overline{\partial}$-closed forms with respect to the metric $\tau_t$ and the fixed K\"ahler metric $\omega.$ Then $P^t(\alpha_t)$ is also Lipschitz, with a Lipschitz constant depending only on that of $\alpha$ and the Lipschitz constant of $\tau_t$ with respect to $t.$
\end{lemma}
   
\noindent In other words, this lemma gives that the integral
$$\int_{X\times\mathbb{D}}\sqrt{-1}dt\wedge d\overline{t}\wedge P^t_\perp(\mu_t^\nu)\wedge\overline{W}e^{-\tau_t^\nu}$$
goes to zero.  Recall that 
$$\partial^{\tau_t^\nu}v_t^\nu=\dot{\tau}_t^\nu\wedge u+h_t^\nu$$
where each $h_t$ is holomorphic on $X$ for each $t$ fixed. Since $\tau_t^\nu$ are non-singular, $v_t^\nu$ depend smoothly on $t.$ Differentiating with respect to $\overline{t}$ yields
$$\partial^{\tau_t^\nu}\frac{\partial v_t^\nu}{\partial\overline{t}}=\left[\frac{\partial^2\tau_t^\nu}{\partial t\overline{\partial}t}-\partial_X(\frac{\partial\tau_t^\nu}{\partial\overline{t}})(V_t)\right]\wedge u+\frac{\partial h_t^\nu}{\partial\overline{t}}.$$
We obtain therefore 
$$P^t_\perp(\mu_t^\nu)=\partial^{\tau_t^\nu}\left(\frac{\partial v_t^\nu}{\partial\overline{t}}\right)$$
since the left hand side is automatically orthogonal to holomorphic forms.
Hence
$$\int_{X\times\mathbb{D}}\sqrt{-1}dt\wedge d\overline{t}\wedge\frac{\partial v_t^\nu}{\partial\overline{t}}\wedge\overline{\overline{\partial}_XW}e^{-\tau_t^\nu}=\int_{X\times\mathbb{D}}\sqrt{-1}dt\wedge d\overline{t}\wedge P^t_\perp(\mu_t^\nu)\wedge\overline{W}e^{-\tau_t^\nu}$$
vanishes as $\nu$ tends to infinity.\\
\medskip

\noindent Now let $\gamma$ be a form of bidegree $(n,1)$ on $X\times\mathbb{D}$ that does not contain any differential $dt.$ Assume that it is Lipschitz with respect to $t$ and decompose it into two parts, $\overline{\partial}_XW,$ which is $\overline{\partial}_X$-exact and the other part orthogonal to $\overline{\partial}_X$-exact forms. We make such orthogonal decomposition for each $t$ separately, by Lemma \ref{le2} each term in the decomposition is still Lipschitz in $t$ and uniformly in $\nu.$ Since $v_t^\nu\wedge\omega^0$ is $\overline{\partial}_X$-closed, it is also the case for $\partial v_t^\nu/\partial\overline{t}.$ We can deduce that $\partial v_t^\nu/\partial\overline{t}$ is $\overline{\partial}$-exact, which is due to Nadel's vanishing theorem \cite{Nad}. Indeed, since 
$$L:=-K_X=c_1(X)-\{\alpha\}+\{\alpha\}$$
and $c_1(X)-\{\alpha\}$ is a K\"ahler class, we can find a K\"ahler form $T^0$ in $c_1(X)-\{\alpha\}$ and a closed positive current $T^1$ such that 
$$T^0+T^1\geq\varepsilon_1\omega^0$$
for some $\varepsilon_1>0.$ The fact that $e^{-2(1+\varepsilon_0)\psi}\in L_{\text{loc}}^1(X)$ implies the multiplier ideal sheaf associated to the local potential of $T^0+T^1$ is trivial. Therefore, by Nadel's vanishing theorem, 
$$H^1(X)=H^1(X, K_X+L)=0.$$
In the sequel, we obtain
$$\int_{X\times\mathbb{D}}\sqrt{-1}dt\wedge d\overline{t}\wedge\frac{\partial v_t^\nu}{\partial\overline{t}} \wedge\overline{\gamma}e^{-\tau_t^\nu} =\int_{X\times\mathbb{D}}\sqrt{-1}dt\wedge d\overline{t}\wedge\frac{\partial v_t^\nu}{\partial\overline{t}}\wedge\overline{\overline{\partial}_XW}e^{-\tau_t^\nu}.$$
Hence
$$\int_{X\times\mathbb{D}}dt\wedge v_t\wedge\overline{\partial_t^{\tau_t^\nu}\gamma}e^{-\tau_t^\nu}$$
goes to zero. Again by Lemma \ref{le1} we can pass to the limit and get 
\begin{equation}\label{eq:21}
\int_{X\times\mathbb{D}}dt\wedge v_t\wedge\overline{\partial_t^{\tau_t}\gamma}e^{-\tau_t}=0
\end{equation}
under the assumption that $\gamma$ is of compact support and Lipschitz in $t.$ We almost show that 
$$\overline{\partial}_tv_t=0$$
except that $\tau_t$ is not smooth. By replacing $\gamma$ by $e^{\tau_t-\tau_t^\prime}\gamma,$ where $\tau_t^\prime$ is another potential on $L.$ If (\ref{eq:21}) holds for some $\tau_t,$ Lipshitz in $t,$ it holds for any such potential. So we can replace $\tau_t$ in (\ref{eq:21}) by some other smooth potential. It shows that $v_t$ is holomorphic in $t$ and therefore $v_t$ is holomorphic in $X\times\mathbb{D}.$\\

\medskip
\noindent Finally, since by definition
$$-v_t=V_t\rfloor u,$$
$V_t$ is holomorphic on $X\times\mathbb{D}.$ Further calculation yields 
$$\overline{\partial}\dot{\tau}_t\wedge u=\partial\overline{\partial}\tau_t\wedge v_t=-\partial\overline{\partial}\tau_t\wedge(V_t\rfloor u)=(V_t\rfloor\partial\overline{\partial}\tau_t)\wedge u.$$
Since $u$ never vanishes, it implies that
$$V_t\rfloor\sqrt{-1}\partial\overline{\partial}\tau_t=\sqrt{-1}\overline{\partial}\dot{\tau}_t.$$   
The Lie derivative of $\partial\overline{\partial}\tau_t$ along $V_t$ is 
\begin{equation}\label{eq:22}
L_{V_t}\partial\overline{\partial}\tau_t=\partial V_t\rfloor\partial\overline{\partial}\tau_t=\partial\overline{\partial}\dot{\tau}_t=\frac{\partial}{\partial t}\partial\overline{\partial}\tau_t.
\end{equation}
Define a holomorphic vector field $\mathcal{V}$ on $X\times\omega$ by
$$\mathcal{V}:=V_t-\frac{\partial}{\partial t}.$$
Let $\eta$ be the form $\partial\overline{\partial}_X\tau_t$ on $\mathcal{X}.$ The formula (\ref{eq:22}) says that
$$L_{\mathcal{V}}\eta=0$$
on $X.$ It follows that $\eta$ is invariant under the flow of $\mathcal{V}$ so $\partial\overline{\partial}\tau_t$ moves by the flow of a holomorphic family of automorphisms of $X.$ In other words, we showed that 
$$F^*(\omega_{\varphi_1}+\Theta)=\omega_{\varphi_0}+\Theta.$$
The equality 
$$\rm{Ric}({\omega_{\varphi_1}})=\omega_{\varphi_1}+\Theta$$
implies that
\begin{equation}\label{eq:23} 
F^*(\rm{Ric}({\omega_{\varphi_1}}))=F^*(\omega_{\varphi_1}+\Theta)=\omega_{\varphi_0}+\Theta=\rm{Ric}({\omega_{\varphi_0}}).
\end{equation}
Since $F^*(\rm{Ric}({\omega_{\varphi_1}}))=\rm{Ric}(F^*(\omega_{\varphi_1})),$ by (\ref{eq:23}), we have  
\begin{equation}\label{eq:24}
\rm{Ric}(F^*(\omega_{\varphi_1}))=\rm{Ric}({\omega_{\varphi_0}}).
\end{equation}
That is to say, there exists a bounded psh function $\zeta$ such that
\begin{equation}\label{eq:25}
F^*(\omega_{\varphi_1})=\omega_{\varphi_0}+\sqrt{-1}\partial\overline{\partial}\zeta.
\end{equation}
(\ref{eq:24}) and (\ref{eq:25}) then induce that 
$$(\omega_{\varphi_0}+\sqrt{-1}\partial\overline{\partial}\zeta)^n=\omega_{\varphi_0}^n$$
and hence
\begin{equation}\label{eq:26}
(\omega+\sqrt{-1}\partial\overline{\partial}(\varphi_0+\zeta))^n=(\omega+\sqrt{-1}\partial\overline{\partial}\varphi_0)^n.
\end{equation}
The following theorem implies that $\zeta=0$ (cf. \cite{EGZ}).

\begin{theorem}
Let $X$ be a compact K\"ahler manifold, $\omega$ a semi positive $(1,1)$-form such that $\int_X\omega^n>0$ and $0\leq f\in L^p(X, \omega^n), p>1,$ a density such that $\int_Xf\phi^n=\int_X\phi^n.$ Then there is a unique bounded function $\varphi$ on $X$ such that $\phi+\sqrt{-1}\partial\overline{\partial}\varphi\geq0$ and 
$$(\omega+\sqrt{-1}\partial\overline{\partial}\varphi)^n=f\phi^n\quad\quad\text{with}\quad\quad\sup_{X}\varphi=-1.$$ 
\end{theorem}

\noindent Therefore (\ref{eq:25}) yields
$$F^*(\rm{Ric}({\omega_{\varphi_1}}))=\rm{Ric}({\omega_{\varphi_0}})$$
and hence $$F^*(\Theta)=\Theta.$$\qed

\noindent \textit {IECN, Universit\'e de Lorraine. B.P. 70239, 54506 Vandoeuvre-l\`es-Nancy Cedex, France}\\
\textit{E-mail adress: ceciliayi$\_$20@hotmail.com}

\end{document}